%%%%%%%%%%%%%%%%%%%%%%%%%%%%%%%%%%%%%%%%%%%%%%%%%%%%%%%%%%%%%%%%%%%%%

%%%
%%%           AMS-TeX 2.1
%%%
%%%
%%%%%%%%%%%%%%%%%%%%%%%%%%%%%%%%%%%%%%%%%%%%%%%%%%%%%%%%%%%%%%%%%%

%%%%%%%%%%%%%%%%%%%%%%%%%%%%%%%%%%%%%%%%%%%%%%%%%%%%%%%%%%%%%%%%%%%%%
\documentstyle{amsppt}
\magnification=\magstep1
\NoRunningHeads

\vsize=7.4in

%%%%%%%%%%%%%%%%%%%%%%%%%%% MACROS   %%%%%%%%%%%%%%%%%%%%%%%%%%%%%%%

 % triplenorm of #1
%%%%%%%%%%%%%%%%%%%%%%%%%%%%%%%%%%%%%%%%%%%%%%%%%%%%%%%%%%%%%%%%%%%%%

\topmatter
\title
Complex interpolation and complementably minimal spaces
\endtitle

\author
P.G. Casazza,
N.J. Kalton, Denka Kutzarova and
M. Masty\l o
\endauthor

\address
Department of Mathematics,
University of Missouri,
Columbia, MO  65211, U.S.A.
\endaddress
\email (Casazza) pete\@casazza.cs.missouri.edu
(Kalton) mathnjk\@mizzou1.bitnet \endemail
\address Permanent address of D. Kutzarova:
Bulgarian Academy of Sciences,
1090 Sofia,
Bulgaria
\endaddress
\email denka\@bgearn.bitnet \endemail
\address Permanent address of M. Masty\l o:
Faculty of Mathematics and Informatics,
\newline
A. Mickiewicz University,
Matejki 48/49, 60-769,
Pozna\'n,
Poland
\endaddress
\email mastylo\@plpuam11.bitnet \endemail

\thanks P. Casazza and N. Kalton were supported by
NSF grant DMS-9201357.
D. Kutzarova was supported by a grant from the Bulgarian
Ministry of Education and Science under contract MM 213-92.
M. Masty\l o was supported by the KBN, grant 2-1051-91-01.
 \endthanks
\abstract

We construct a class of super-reflexive complementably minimal
spaces, and study uniformly convex distortions of the norm on Hilbert
space by using methods of complex interpolation.

\endabstract

\endtopmatter

\document \baselineskip=14pt
\heading{1. Introduction}
\endheading

A Banach space $X$ is called {\it (complementably) minimal} if every
infinite-dimensional closed subspace $E$ contains a
(complemented) subspace  isomorphic to $X.$  These notions were
introduced by Pe\l czy\'nski \cite{17} and Rosenthal \cite{19}. Any
minimal
space
must be separable and it is classical that the spaces $\ell_p$ for $1\le
p<\infty$ and $c_0$ are complementably minimal. The space $T^*$ (the dual
of Tsirelson space) provides another example of a minimal but not
complementably minimal space (\cite{2}, \cite{3}).  Recently
Schlumprecht
(\cite{20},
\cite{21})
constructed the first example of a complementably minimal space other
than the classical spaces $\ell_p$ and $c_0$ and this was a launching
point for a number of remarkable developments in Banach space theory
(\cite{8},\cite{15},\cite{16}).

The space constructed by Schlumprecht is reflexive (see
\cite{20} and Proposition 2 below) but fails to
be super-reflexive since it contains $\ell_{\infty}^n$'s uniformly.
Our main aim in this note is to show how interpolation methods can be
used to extend Schlumprecht's construction and thereby introduce a  class
of complementably minimal super-reflexive spaces.  We also show that
interpolation can be used to tighten the known results on distortions of
the norm in Hilbert space (\cite{16}); precisely one can require the
distorted
norm to satisfy good uniform convexity and uniform smoothness conditions.

Our arguments depend heavily on the ideas of complex interpolation of
Banach spaces first introduced by Calder\'on in 1964 (\cite{1}).  For
another application of such ideas to problems of this nature see Daher
\cite{5}.

\heading{2.  Some remarks on minimal and complementably minimal spaces}
\endheading

 \proclaim{Proposition 1}Every minimal Banach space is isomorphic to a
subspace of a minimal space with an unconditional basis.\endproclaim

\demo{Proof}If $X$ is minimal and contains an unconditional basic
sequence then it is clear that $X$ embeds into a  minimal space with
unconditional basis.  We argue that $X$ must contain an unconditional
basic sequence.  For otherwise by \cite{7} $X$ contains an hereditarily
indecomposable subspace $Y$.  Since $Y$ is not isomorphic to any proper
subspace of itself (\cite{8} Corollary 19 and Theorem 21) this is a
contradiction.\enddemo

\proclaim{Proposition 2}(1) A minimal space is either reflexive or
isomorphic to a subspace of $c_0$ or $\ell_1.$ \newline (2) A
complementably minimal space is either
$c_0$,
$\ell_1$ or is reflexive.\endproclaim

\demo{Proof}(1) follows immediately from the preceding
Proposition, since any space with an unconditional basis contains $c_0$
or $\ell_1$ or is reflexive.  (2) follows quickly from (1).
\enddemo

We remark that it is unknown if every complementably minimal space
is prime, or whether complementable minimality passes to dual spaces in
general.  Note that the minimal space $T^*$ has a non-minimal dual
(\cite{3}).

\heading{3. Some classes of complementably minimal spaces}
\endheading

We let $c_{00}$ be the space of all finitely non-zero sequences.
If $E_1$ and $E_2$ are finite intervals of natural numbers  we write
$E_1<E_2$ to mean $\max E_1<\min E_2.$  If $x\in c_{00}$ and $E$ is a
subset of $\bold N$ we write $Ex=x\chi_E.$  We will also need the concept
of a {\it block subspace}: this is a subspace of $c_{00}$ generated by a
sequence $(u_n)$ whose supports $E_n$ satisfy $E_1<E_2<\cdots.$

We will consider spaces $X$ determined by lattice norms $\|.\|_X$ on
$c_{00}.$  We will then let $X$ be the space of all sequences $x$ so that
$\|x\|_X=\sup\|(x_1,\ldots,x_n,0,\ldots)\|_X<\infty.$  We abbreviate
$\|x\|_{\ell_p}$ to $\|x\|_p$ if $1\le p\le \infty.$  If
$X$ and
$Y$ are
two such spaces and $0<\theta<1$ we define $X^{1-\theta}Y^{\theta}$ to be
the space $Z$ defined by
$\|z\|_Z=\inf\{\max(\|x\|_X,\|y\|_Y):\ |z|=|x|^{1-\theta}|y|^{\theta}\}.$
When working over the complex scalars, if either $X$ or $Y$ is
separable  then
$Z$
coincides with the usual complex interpolation space $[X,Y]_{\theta}$
(see \cite{1}).  It will, however, be easily seen that our results
apply also in the real case.

We let $\Cal G$ be the class of increasing functions $f:[1,\infty)\to
[1,\infty)$
so that:\newline
(1) $f(1)=1$ and $f(x)<x$ if $x> 1.$ \newline
(2) $x/f(x)$ is concave.\newline
(3) $f$ is submultiplicative, i.e. $f(xy)\le f(x)f(y).$ \newline
Suppose further $  1\le p<r\le \infty,$ and $f\in\Cal G.$ We define $\Cal
X(p,r;f)$ to be collection
of all sequence spaces $X$ so that:\newline (4) $\|x\|_r\le
\|x\|_X\le
\|x\|_p$ for all $x\in c_{00}.$ \newline
(5) $X$ is $p$-convex and $r$-concave (with constants
one).\newline
(6) If $E_1<E_2<\cdots<E_n$ are intervals in $\bold N$ then
$$ \|x\|_X \ge
\frac{1}{f(n)^{1/p-1/r}}(\sum_{i=1}^n\|E_ix\|_X^p)^{1/p}.$$
Then $\Cal X(p,r;f)$ is non-empty.  Furthermore if $r<\infty$ then any
$X\in\Cal X(p,r;f)$ is separable by $r$-concavity.  If $r=\infty$ the
same conclusion can be obtained from the fact that $f(n)=o(n)$
as
$n\to
\infty,$ since $X$ cannot then contain $c_0.$  When
$p=1$ and
$r=\infty$
 it is clear that there is a unique space, which may be
constructed by an inductive procedure (\cite{20},\cite{8})
$S=S(f)$
(Schlumprecht $f$-space) satisfying a minimality condition:
 $\|x\|_S\le \|x\|_X$ for $X\in\Cal X(1,\infty;f).$
Furthermore as shown by Schlumprecht (\cite{20}) in this space
$\|\sum_{i=1}^n e_i\|_S = n/f(n)$ where $e_i$ are the basis vectors.

\proclaim{Proposition 3}For any $f\in\Cal G$ and $1\le p<r\le \infty,$
the space
$S_{p,r}=\ell_t^{1-\theta}S^{\theta}$ where
$\theta=\frac1p-\frac1r$ and
$t=(1-\theta)r,$ is the
unique space in $\Cal X(p,q;f)$ satisfying $\|x\|_{S_{p,r}}\le \|x\|_X$
for all
$X\in\Cal X(p,r;f).$  Furthermore
$\|\sum_{i=1}^ne_i\|_{S_{p,r}}=n^{1/p}f(n)^{1/r-1/p}.$\endproclaim

\demo{Proof}For the case $r=\infty$ this follows by elementary
convexification from the case $p=1,$ since the space $S_{p,\infty}$
coincides with the $p$-convexification of $S.$  If
$r<\infty$ it will follow easily by
convexification or concavification from the case when $p<2$ and $r=q$ the
conjugate
index of $p.$  We therefore suppose $r=q$.  First we prove that the
space $S_{p,q}$ is in the class $\Cal X(p,q;f).$
Indeed the only property to be verified is (3), and this is standard.  If
$0\le x\in c_{00}$ and
$E_1<E_2<\cdots<E_n$ then we can write $x=u^{1-\theta}v^{\theta}$
where $0\le u,v\in c_{00}$ and $\|u\|_2=\|v\|_S=\|x\|_{S_{p,q}}.$  Thus
$$
\align
(\sum_{j=1}^n\|E_jx\|_{S_{p,q}}^p)^{1/p} &\le
(\sum_{j=1}^n\|E_ju\|_2^{p(1-\theta)}\|E_jv\|_S^{p\theta})^{1/p}\\
&\le
(\sum_{j=1}^n\|E_ju\|_2^2)^{1/q}(\sum_{j=1}^n\|E_jv\|_S)^{1/p-1/q}\\
&\le f(n)^{1/p-1/q}\|u\|_2^{2/q}\|v\|_S^{1/p-1/q}=
f(n)^{1/p-1/q}\|x\|_{S_{p,q}}.
\endalign
$$

  Conversely suppose
$X\in\Cal X(p,q;f)$. Then
by Pisier's extrapolation theorem (\cite{18} or \cite{9}) there is a
sequence
space $Y$ so that $X=Y^{\theta}\ell_2^{1-\theta}$ where $\theta
=\frac1p-\frac1q=2(\frac1p-\frac12).$  We show that $Y\in\Cal
X(1,\infty;f).$  Clearly $\|y\|_{\infty}\le \|y\|_Y\le \|y\|_1$ for all
$y\in c_{00}.$  Now suppose $0\le y\in c_{00}$ and $E_1<E_2<\cdots<E_n$
are disjoint intervals; let $y_i=E_iy$.  Pick $y_i^*\in c_{00}$ supported
on
$E_i$ with
$\|y_i^*\|_{Y^*}=1$ and
$\langle y_i,y_i^*\rangle=\|y_i\|_Y.$  Let
$v_i=y_i^{\frac{1+\theta}2}(y_i^*)^{\frac{1-\theta}2}=y_i^{\theta}
w_i^{1-\theta}$ where $w_i=(y_iy_i^*)^{1/2}.$  Then
$\|w_i\|_2=\|y_i\|^{1/2}_Y.$

Now $X=Y^{\frac{1+\theta}2}(Y^*)^{\frac{1-\theta}2}$ (this
 follows for example
from Lozanovskii's theorem (\cite{12}) and the Re-iteration theorem
(\cite{1} and \cite{4})). Also by
the duality theorem $X^*=Y^{\frac{1-\theta}2}(Y^*)^{\frac{1+\theta}2}.$
Hence $\|v_i^*\|_{X^*}\le \|y_i\|_Y^{\frac{1-\theta}2}$ where
$v_i^*=y_i^{\frac{1-\theta}2}(y_i^*)^{\frac{1+\theta}2}.$  Hence
$\|v_i\|_X\ge \|y_i\|_Y^{\frac{1+\theta}2}.$

We conclude that
$$ \|\sum_{i=1}^nv_i\|_X \ge \frac{1}{f(n)^{\theta}}(\sum_{i =1}^n
\|y_i\|_Y^{p(1+\theta)/2})^{1/p}.$$
This implies
$$ \|\sum_{i=1}^nv_i\|_X \ge
\frac{1}{f(n)^{\theta}}(\sum_{i=1}^n\|y_i\|_Y)^{(1+\theta)/2}.$$

On the other hand
$$ \|\sum_{i=1}^nv_i\|_X \le \|y\|_Y^{\theta} \|\sum w_i\|_2^{1-\theta}$$
or
$$ \|\sum_{i=1}^nv_i\|_X\le \|y\|_Y^{\theta} (\sum_{i=1}^n \|y_i\|_Y)^{(1
-\theta)/2}.$$
Combining these inequalities gives that
$$ \|y\|_Y \ge \frac1{f(n)}\sum_{i=1}^n\|y_i\|_Y$$
so that $Y\in\Cal X(1,\infty;f).$  Thus $Y\subset S$ with norm-one
inclusion and similarly $ Y^{\theta}\ell_2^{1-\theta}=X\subset S_{p,q}$
with norm one inclusion either by interpolation or simple
calculation of norms.  Finally we note that
$$\|\sum_{i=1}^ne_i\|_{S_{p,r}}\le
\|\sum_{i=1}^ne_i\|_2^{1-\theta}\|\sum_{i=1}^ne_i\|_S^{\theta}=
n^{1/p}f(n)^{1/p-1/q}.$$
The other inequality follows from (6) immediately since the basis vectors
have norm one.\enddemo

\demo{Remark}The last condition shows that $S_{p,r}$ does not coincide
with $\ell_p$ as a sequence space.  It further follows that since the
basis $(e_n)$ of $S_{p,r}$ is subsymmetric that it has no subsequence
equivalent to the $\ell_p-$basis and therefore $S_{p,r}$ cannot be even
isomorphic to $\ell_p.$
\enddemo

The following remarkable result is due to Schlumprecht \cite{21}:

\proclaim{Theorem 4}If $f(x)=\log_2(x+1)$ then $S=S(f)$ is
complementably minimal.\endproclaim

We now prove a simple extension of a technique used by both Schlumprecht
(see \cite{20} and \cite{21}) and Gowers and Maurey (\cite{8}).

\proclaim{Proposition 5}Suppose $f\in \Cal G$ and that
 $\lim_{x\to \infty}f(x)x^{-a}=0$ if $a>0.$
Suppose $1\le p< \infty,$ and suppose $\frac1p+\frac1q=1$ (with
appropriate interpration when $p=1).$ Suppose $X\in \Cal X(p,\infty;f)$
Then:\newline
(1) If $n\in\bold N$ and $\epsilon>0$ and $W$ is a block subspace of
$c_{00}$ there is a block basic sequence
$(u_1,u_2,\ldots,u_n)$ in $W$ so that $\|u_i\|_X= 1$ for $1\le i\le
n$ and $\|u_1+\cdots+u_n\|_X\ge n^{1/p}-\epsilon.$\newline
(2) If $n\in\bold N$ and $\epsilon>0$ and $W$ is a block subspace of
$c_{00}$ there is a block basic sequence
$(u^*_1,u^*_2,\ldots,u^*_n)$ in $W$ so that $\|u^*_i\|_{X^*}= 1$ for
$1\le
i\le n$ and $\|u^*_1+\cdots+u_n^*\|_{X^*}\le n^{1/q}+\epsilon$.\newline
\endproclaim

\demo{Proof}(1) can be obtained immediately from Lemma 3 of Gowers-Maurey
(\cite {8}) by convexification.  The proof of (2) is similar.  Let
$\beta_n=\beta_n(W)$ be the least constant so that for every $k,\epsilon$
there
exists a normalized block basic sequence $(u_1^*,u_2^*,\ldots,u_n^*)$ in
$W$ with $k<$ supp $u_1^*$ and $\|\sum_{i=1}^nu_i^*\|_{X^*}\le
\beta_n+\epsilon.$  Then it is easy to see that $\beta_{mn}\ge
\beta_m\beta_n$ for $m,n.$  Also by $q$-concavity of $X^*$ we have
$\beta_n\ge n^{1/q}$ for all $n.$  On the other hand one can verify
easily by duality that
for all block basic sequences $(u_1^*,\ldots,u_n^*)$ we have
$\|\sum_{i=1}^nu_i^*\|_{X^*}\le
f(n)(\sum_{i=1}^n\|u_i^*\|_{X^*}^q)^{1/q}$ and so $\beta_n\le
f(n)n^{1/q}$.  Hence $\beta_n\le (\beta_{n^k})^{1/k}\le
(f(n^k))^{1/k}n^{1/q}$ and so $\beta_n=n^{1/q}.$\enddemo

The key to Schlumprecht's argument for Theorem 4 is the following, which
combines his Lemma 2 and Theorem 3.  (Note that if $(u_i)_{i=1}^n$ is a
normalized block basic sequence with $\|u_1+\cdots+u_n\|>n-\epsilon$ then
$(u_i)_{i=1}^n$ is $(1-\epsilon)^{-1}$ equivalent to the $\ell_1^n$
basis.)

\proclaim{Proposition 6}(Schlumprecht \cite{21}) Suppose
$f(x)=\log_2(x+1)$
and that $(u_n)_{n=1}^{\infty}$ is a normalized block basic sequence in
$S=S(f).$
 Let $v_n=2^{-n}\sum_{i=2^n+1}^{2^{n+1}}u_i$ and suppose
$\lim_{n\to\infty}2^n(1-\|v_n\|_S)=0.$  Then $(v_n)$ has
a subsequence $(w_n)$ so that $(w_n)$ is equivalent to the unit vector
basis of $S(f).$\endproclaim

\demo{Remark}The conclusions of Propositions 5 and 6 are all we require
for our main result.  Thus Theorem 8 below will hold for any $f\in\Cal F$
for which these Propositions hold; clearly this is a much wider class
than just the singleton $f(x)=\log_2(x+1)$ but it has not been precisely
determined to date.\enddemo

\proclaim{Proposition 7}Suppose $f(x)=\log_2(x+1)$ and $1\le
p<r\le
\infty$. Suppose $(u_n)$ is a normalized block basic sequence in
$S_{p,r}=S_{p,r}(f).$
Let $v_n=2^{-n/p}\sum_{i=2^n+1}^{2^{n+1}}u_i$ and suppose that
$\lim_{n\to\infty}2^n(1-\|v_n\|_{S_{p,r}})=0.$   Then $(v_n)$ has a
subsequence
$(w_n)$
which is equivalent to the unit vector basis of $S_{p,r}.$\endproclaim

\demo{Proof}If $r=\infty$ this follows immediately from the fact that
$S_{p,\infty}$ is the $p$-convexification of $S.$  If $r<\infty$ we
can assume each $u_i\ge 0.$  Let
$\epsilon_n=2^{n}(1-\|v_n\|_{S_{p,r}}),$ so that
$\lim_{n\to\infty}\epsilon_n=0.$ Note that each
$u_i=x_i^{\theta}y_i^{1-\theta}$ where $(x_i)$ is a normalized positive
block basic sequence in $S$ and $(y_i)$ is a normalized block basic
sequence in $\ell_t$ as in Proposition 3.  Then $$\|v_n\|_{S_{p,r}} \le
2^{-n/p}\|\sum_{i=2^n+1}^{2^{n+1}}x_i\|_S^{\theta}\|\sum_{i=2^n+1}^
{2^{n+1}}y_i\|_t^{1-\theta}$$
so that
$$ 1-2^{-n}\epsilon_n \le
2^{n((1-\theta)/t-1/p)}\|\sum_{i=2^n+1}^{2^{n+1}}x_i\|_S^{\theta}$$
and this simplifies to
$$ \|\sum_{i=2^n+1}^{2^{n+1}}x_i\|_S\ge
(2^{n\theta}-2^{n(\theta-1)}\epsilon_n)^{1/\theta}\ge
2^n-\theta^{-1}\epsilon_n.$$

It thus follows from Proposition 6 that a suitable subsequence
$(x_{k(n)})$ is equivalent
to the unit vector basis of $S$ and so for some constant $K$ and any
finitely nonzero sequence $(d_n)$ we have
$$ \|\sum_{n=1}^{\infty}d_nx_{k(n)}\|_S \le K\|d\|_S.$$
Now if $0\le a,b\in c_{00}$ then
$$ \|\sum_{n=1}^{\infty}a_n^{\theta}b_n^{1-\theta}v_{k(n)}\|_{S_{p,r}}
\le \|\sum a_n x_{k(n)}\|_S^{\theta} \|\sum b_n y_{k(n)}\|_t^{1-\theta}
$$
Hence
$$ \|\sum_{n=1}^{\infty}a_n^{\theta}b_n^{1-\theta}v_{k(n)}\|_{S_{p,r}}
\le K^{\theta}\|a\|_S^{\theta}\|b\|_t^{1-\theta}$$
and so for any $d\in c_{00}$
$$ \|\sum_{n=1}^{\infty} d_n v_{k(n)}\|_{S_{p,r}} \le
K^{\theta}\|d\|_{S_{p,r}}.$$

However the minimality property of the norm on $S_{p,r}$ clearly implies
that
$$ \|d\|_{S_{p,r}} \le \|\sum_{n=1}^{\infty}d_n
v_{k(n)}\|v_{k(n)}\|_{S_{p,r}}^{-1}\|_{S_{p,r}} \le
C_0\|\sum_{n=1}^{\infty} d_n v_{k(n)}\|_{S_{p,r}},$$
where $C_0=\sup_n \|v_{k(n)}\|_{S_{p,r}}^{-1}.$
Taking $w_n=v_{k(n)}$ the result is proved.\enddemo

\proclaim{Theorem 8}Suppose $f(x)=\log_2(x+1)$ and that $1\le
p<r\le
\infty.$
Then the spaces $S_{p,r}(f)$ and $S_{p,r}(f)^*$ are complementably
minimal.\endproclaim

\demo{Proof} By combining Proposition 5 and Proposition 7 we see that any
block subspace of $S_{p,r}$ contains a normalized block basic sequence
$(w_n)$ equivalent to the unit vector basis of $S_{p,r}.$  Let $(g_n)$ be
any bounded block basic sequence in $S_{p,r}^*$ so that $E_n=$ supp $g_n$
does not meet supp $w_m$ when $m\neq n$ and $\langle w_n,g_n\rangle =1.$
Then we claim that $Px=\sum_{n=1}^{\infty}\langle x,g_n\rangle w_n$
defines a projection onto $[w_n].$  In fact if $x\in c_{00}$ then for
a suitable constant $C,$ and letting $M=\sup\|g_n\|_{S_{p,r}^*},$
$$
\align
\|Px\|_{S_{p,r}} &\le C \|(\langle x,g_n\rangle)\|_{S_{p,r}}\\
       &\le CM \|(\|E_nx\|)\|_{S_{p,r}}\\
       &\le CM\|x\|_{S_{p,r}}
\endalign
$$
and so $[w_n]$ is complemented. This shows that $S_{p,r}$ is
complementably minimal.

In $S_{p,r}^*$ we argue that any block subspace contains a
normalized block basic sequence $(u_n^*)$ so that
$\|2^{-n/q}\sum_{i=2^n+1}^{2^{n+1}}u_i^*\|_{S_{p,r}^*}\le
(1+\frac1n2^{-n})
$ for all
$n.$ Choose
a normalized block basic sequence $(u_n)$ in $S_{p,r}$ so that supp $u_n$
is contained in supp $u_n^*$ and $\langle u_n,u_n^*\rangle =1.$  Then
$\|2^{-n/p}\sum_{i=2^n+1}^{2^{n+1}}u_i\|\ge 1-\frac1n2^{-n}$.  It
follows by the
argument of the first part that we can select a sequence $k(n)\to \infty$
so that if $w_n = 2^{-k(n)/p}\sum_{i=2^{k(n)}+1}^{2^{k(n)+1}}u_i$ then
$(w_n)$ is equivalent to the unit vector basis of $S_{p,r}$ and
complemented by the projection
$Px= \sum_{n=1}^{\infty}\langle x,w_n^*\rangle w_n$ where
$w_n^*=2^{-k(n)/q}\sum_{i=2^{k(n)+1}}^{2^{k(n)+1}}u_i^*.$  Thus $[w_n^*]$
is a complemented subspace of $S_{p,r}^*$ isomorphic to
$S_{p,r}^*.$\enddemo

\demo{Remark}Of course in the cases when $1<p<r<\infty$ the space
$S_{p,r}$ is super-reflexive (in fact the given norm is uniformly
convex).
Note also that it is trivial that none of these spaces can contain a copy
of $c_0$ or $\ell_p$ for any $1\le p<\infty.$  The first such example
was due to Tsirelson (\cite{22}) and the first super-reflexive
example of such a space was given by Figiel and Johnson (\cite{6}).
Note also that
$S^*=S_{1,\infty}^*$ is complementably minimal.
We remark also that the spaces $S_{p,r}$ for $1<p<r<\infty$ are
arbitrarily distortable in view of their minimality and a recent result
of Maurey
\cite{13} (see also \cite{14}).
\enddemo

\heading {4. Distortions of Hilbert space}\endheading

The distortion problem for Hilbert spaces was recently solved by Odell
and Schlumprecht (\cite{15},\cite{16}); for the latest developments see
Milman and Tomczak-Jaegermann \cite{14} and Maurey \cite{13}.
In this brief section we show the following, a strengthening of Theorem
1.2 of \cite{16}.

\proclaim{Theorem 9}Given $1<p<2$ (with conjugate index $q$) there is a
constant
$K=K(p)$ so that
for any constant $C$ there is a Banach space $X$ isomorphic to $\ell_2$
 so that:\newline
(1) $X$ contains no $C$-unconditional basic
sequence.\newline
(2) The modulus of convexity $\delta_X$ satisfies
$\delta_X(\epsilon)\ge K^{-1}\epsilon^q$ for $0\le\epsilon\le 1.$
\newline
(3) The modulus of smoothness $\rho_X$ satisfies
$\rho_X(\tau)\le K\tau^p$ for $0\le \tau\le 1.$
\endproclaim

\demo{Proof}In fact this is an easy modification of Gowers-Maurey
\cite{8} Theorem 2, using the Odell-Schlumprecht theorem on the
existence of asymptotic biorthogonal systems in $\ell_2$ (\cite{15},
\cite{16}).
We shall suppose that all spaces are complex in the following argument.
It is clear that the construction actually yields a space in which the
underlying real space has no good unconditional basic sequences, so that
the complex structure can be ``forgotten'' at the end.

Fix $\theta=2/q$ and choose $r$ an integer so that
$r^{1-\theta}>8C.$
Using the notation of \cite{8} let $(A_n,A_n^*)$ be an asymptotic
biorthogonal system with constant at most $r^{-2}.$  For each $n$ let
$Z_n^*$ be a countable dense subset of $A_n^*$ and let $Z^*=\cup Z_n^*.$
Let $\sigma$ be an injection from the collections of all finite sequences
in $Z^*$ to $\bold N.$  Let $\Gamma_r$ be the
collection of all
$z^*\in \ell_2^*$ of the form $z^*=\sum_{i=1}^rz_i^*$ where $z_1^*\in
Z_1^*$ and $z_{i+1}^*\in Z^*_{\sigma(z_1^*,\ldots,z_i^*)}.$

Define a new equivalent norm $\|.\|_Y$ on $\ell_2$ by
$\|x\|_Y=\max(\|x\|,r\max\{|\langle x,z^*\rangle|:z^*\in
\Gamma_r\}).$
Let $Y=(\ell_2,\|.\|_Y)$ and form the complex interpolation
space
$X=[Y,\ell_2]_{\theta}.$

Note first that $X$ is $\theta-$Hilbertian (\cite{18}) and so
satifies the inequalities
$$ \frac12(\|x+y\|_{X}^p+\|x-y\|_{X}^p) \le
\|x\|_{X}^p+\|y\|_{X}^p$$
and
$$ \frac12(\|x+y\|_{X}^q+\|x-y\|_{X}^q) \ge
\|x\|_{X}^q+\|y\|_{X}^q$$
>From this  it follows easily that one has the appropriate
estimates on the moduli of convexity and smoothness.

It remains to show that $X$ has no $C$-unconditional basic sequence.
To see this let $E$ be any block subspace of $X$.  Arguing as in
\cite{8} one constructs a block basis
$(z_1,\ldots,z_r)$ normalized in $\ell_2$ and $z^*\in\Gamma_r$ so that
$|\langle\sum_{i=1}^rz_i,z^*\rangle| \ge r-1$ and
$\|\sum_{i=1}^r(-1)^iz_i\|_Y\le 4r.$

Now $\|z^*\|_{Y^*}\le r^{-1}$ and $\|z^*\|_2 \le r^{1/2}.$   Hence, by
the duality theorem
$\|z^*\|_{X^*}\le r^{(3\theta-2)/2}$, whence it follows that
$\|\sum_{i=1}^r z_i\|_X \ge \frac12 r^{(4-3\theta)/2}.$  On the other
hand
$\|\sum_{i=1}^r (-1)^iz_i\|_2 \le r^{1/2}$ so that
$\|\sum_{i=1}^r (-1)^iz_i\|_X \le 4^{1-\theta}r^{(2-\theta)/2}$.  It
follows no basic sequence in $X$ has unconditional basis constant better
than $\frac18 r^{1-\theta}>C$.
\enddemo

\demo{Remark} In this theorem one cannot achieve a similar result with
$p=q=2.$  For in this case it follows that $X$ has bounded type two and
cotype two constants (cf. \cite{11} p. 77) and so by
Kwapie\'n's theorem (\cite{10}) has bounded distance to Hilbert space.
\enddemo


\Refs

\ref\no{1}\by A.P. Calder\'on \paper Intermediate spaces and
interpolation, the complex method \jour Studia Math. \vol 24 \yr 1964
\pages 113-190\endref

\ref\no{2}\by P.G. Casazza, W.B. Johnson and L. Tzafriri \paper On
Tsirelson's space \jour Israel J. Math. \vol 47\yr 1974 \pages
191-218\endref

\ref\no {3} \by P.G. Casazza and T.J. Shura \book Tsirelson's space
\bookinfo Springer Lecture Notes 1363 \publ Springer\publaddr New York
\yr 1988\endref

\ref\no{4}\by M. Cwikel \paper Complex interpolation spaces, a discrete
definition and reiteration \jour Indiana Univ. Math. J. \vol 27 \yr 1978
\pages 1005-1009\endref

\ref\no{5}\by M. Daher \paper Hom\'eomorphismes uniformes entre les
sph\`eres unit\'e des espaces d'interpolation \jour C.R. Acad. Sci.
(Paris) \vol 316 \yr 1993 \pages 1051-1054\endref

\ref\no{6}\by T. Figiel and W.B. Johnson\paper A uniformly convex space
which contains no $\ell_p$ \jour Comp. Math. \vol 29\yr 1974 \pages
179-190\endref

\ref\no {7} \by W.T. Gowers \paper A new dichotomy for Banach spaces
\paperinfo preprint
\endref

\ref\no {8} \by W.T. Gowers and B. Maurey \paper The unconditional basic
sequence problem \jour J. Amer. Math. Soc. \vol 6 \yr 1993\pages
851-874\endref


\ref\no{9} \by N.J. Kalton \paper Differentials of complex
interpolation processes for K\"othe function spaces \jour Trans. Amer.
Math. Soc. \vol 333 \yr 1992 \pages 479-529
\endref

\ref\no{10}\by S. Kwapie\'n \paper Isomorphic charcaterization of
inner-product spaces by orthogonal series with vector-valued coefficients
\jour Studia Math. \vol 44 \yr 1972 \pages 583-595  \endref

\ref\no{11}\by J. Lindenstrauss and L. Tzafriri \book Classical Banach
spaces II, Function spaces \publ Springer\publaddr Berlin \yr
1979\endref

\ref\no {12}  \by G.Y. Lozanovskii \paper On some Banach lattices \jour
Siberian Math. J. \vol 10 \yr 1969 \pages 419-430\endref

\ref\no{13} \by B. Maurey \paper A remark about distortion \paperinfo
preprint\endref

\ref\no{14} \by V.D. Milman and N. Tomczak-Jaegermann \paper Asymptotic
$\ell_p$-spaces and bounded distortions \jour Contemp. Math. \vol
144 \yr 1993 \pages 173-195\endref

\ref\no {15} \by E. Odell and T. Schlumprecht \paper The distortion of
Hilbert space \jour Geom. and Funct. Anal. \vol 3\yr 1993 \pages 201-217
\endref

\ref\no {16} \by E. Odell and T. Schlumprecht \paper The distortion
problem \paperinfo preprint\endref

\ref\no{17} \by A. Pe\l czy\'nski \paper Projections in certain Banach
spaces \jour Studia Math. \vol 19 \yr 1960 \pages 209-228\endref

\ref\no {18} \by G. Pisier \paper Some applications of the complex
interpolation method to Banach lattices \jour J. d'Analyse Math. \vol 35
\yr 1979 \pages 264-281\endref


\ref\no{19} \by H.P. Rosenthal \paper On a theorem of Krivine concerning
block finite representability of $\ell_p$ in general Banach spaces \jour
J. Functional Analysis \vol 28 \yr 1978 \pages 197-225\endref

\ref\no{20}\by T. Schlumprecht \paper An arbitrarily distortable
Banach space \jour Israel J. Math. \vol 76 \yr 1991 \pages 81-95\endref

\ref\no {21} \by T. Schlumprecht \paper A complementably minimal space
not containing $c_0$ or $\ell_p$ \paperinfo preprint \endref

\ref\no{22}\by B.S. Tsirelson \paper Not every Banach space contains an
embedding of $\ell_p$ or $c_0$ \jour Funct. Anal. Appl. \vol 8\yr
1974\pages 138-141\endref


\endRefs
\enddocument